\documentclass[reqno, 12pt]{amsart}
\usepackage{amsmath,amsfonts,amsthm,amscd,upref,amstext}
\usepackage{amssymb}

\usepackage{subfigure}
\usepackage[all]{xy}
\usepackage{comment}
\usepackage{xcolor}
\usepackage{fullpage}
\usepackage[colorinlistoftodos]{todonotes}

\usepackage{calrsfs}

\newtheorem{prop}{Proposition}[section]
\newtheorem{thm}[prop]{Theorem}
\newtheorem{cor}[prop]{Corollary}
\newtheorem{ex}[prop]{Example}
\newtheorem{defn}[prop]{Definition}

\newtheorem{conjec}[prop]{Conjecture}
\newtheorem{rem}[prop]{Remark}

\newtheorem{nota}[prop]{Notation}
\newtheorem{fact}[prop]{Fact}

\newcommand\edim{\operatorname{edim}}
\newcommand\mult{\operatorname{mult}}

\newcommand\pd{\operatorname{pd}}

\newcommand{\m}{\mathfrak{m}}

\usepackage{tikz}

\usepackage{tikz-cd}
\usetikzlibrary{matrix}

\usepackage{lmodern} 
\usepackage{amsmath}
\usepackage{amssymb}
\usepackage{graphicx,scalerel}

\usepackage{scalerel,stackengine}
\stackMath

\newcommand\reallywidehat[1]{%
\savestack{\tmpbox}{\stretchto{%
  \scaleto{%
    \scalerel*[\widthof{\ensuremath{#1}}]{\kern-.8pt\bigwedge\kern-.8pt}%
    {\rule[-\textheight/2]{1ex}{\textheight}}
  }{\textheight}%
}{0.85ex}}%
\stackon[1pt]{#1}{\tmpbox}%
}

\numberwithin{equation}{section}

\begin{document}

\title[Structural results of fiber product rings]{On General fiber product rings, Poincar\'e series and their structure}

\date{}


\author{T. H. FREITAS}
\address{Universidade Tecnol\'ogica Federal do Paran\'a, 85053-525, Guarapuava-PR, Brazil}
\email{freitas.thf@gmail.com}


\author{J. A. LIMA}
\address{Universidade Tecnol\'ogica Federal do Paran\'a, 85053-525, Guarapuava-PR, Brazil}
\email{seyalbert@gmail.com}



\begin{abstract}
The present paper deals with the investigation of the structure of general fiber product rings $R\times_TS$, where $R$, $S$ and $T$ are local rings with common residue field. We show that the Poincar\'e series of any $R$-module over the fiber product ring  $R\times_TS$ is bounded by a rational function. In addition,  we give a description of ${\rm depth}(R\times_TS)$, which is an open problem in this theory.  As a biproduct, using the characterization of the Betti numbers over $R\times_TS$ obtained,   we provide certain cases of   the Cohen-Macaulayness of  $R\times_TS$ and, in particular, we show that  $R\times_TS$ is always non-regular.    Some positive answers for the Buchsbaum-Eisenbud-Horrocks and Total rank conjectures over $R\times_TS$ are also established.        
\end{abstract}

\maketitle

\section{Introduction}  
Let $(R,\mathfrak{m}_R,k)$, $(S,\mathfrak{m}_S,k)$ and $(T,\mathfrak{m}_T,k)$ be
 commutative Noetherian local rings,  and let $R \stackrel{\pi_R}\twoheadrightarrow T \stackrel{\pi_S}\twoheadleftarrow S$ be
  surjective homomorphisms of rings. Recently, the study of the fiber product rings $R\times_TS$ has been an effective topic of investigation by several authors in different  subjects, for instance, structural results \cite{AAM, shirogoto, grow, nstv, ramati}, homological problems \cite{paperbetti, paper, nasseh, hugh, moore, nstv, NT}, and their geometric aspects 
 \cite{paperanalytic, paperformal}. Actually, most of these works consider the case $T=k$, and this  provides a good understanding of the structure of $R\times_kS$. When $T\neq k$, few results are known and this motivate us the investigation showed here.  

 The main focus of this paper is, roughly speaking, the study of the structure of the fiber product ring $R\times_TS$, precisely, when it is Cohen-Macaulay, complete  intersection, hypersurface and a regular ring. In fact, the structure of the fiber product ring can change severely, depending  of the structure of the rings $R$, $S$ and $T$. For instance, if $R$ and $S$ are $d$-dimensional regular rings, with $d\geq 2$, then $R\times_kS$ is also $d$-dimensional, but is not Cohen-Macaulay (see Fact \ref{dimdepth}). Also,   $R\times_kS$ is Gorenstein if and only if it is a $1$-dimensional hypersurface \cite[Remark 3.1]{grow} (or \cite[Corollary 2.7]{nstv}).

 For this purpose, we fix the investigation in three key topics:   Poincar\'e series,  Betti numbers and the depth of the fiber product ring. First, we show that the Poincar\'e series $P^{R\times_TS}_M(t)$ of any $R$-module $M$ over the fiber product ring $R\times_TS$ is bounded by a rational function which depends on $P^R_M(t)$, $P^R_T(t)$ and $P^S_T(t)$. This main result is given below (Theorem \ref{thm1}): 
 \begin{thm} \label{thm10} Let  $R\times_TS$ be a  fiber product ring. If $M$ is an $R$-module, then
$$P^{R\times_TS}_M(t) \succeq \frac{P^R_M(t)P^S_T(t)}
{P^R_T(t)+P^S_T(t)-P^R_T(t)P^S_T(t)}.$$
\end{thm}

As a biproduct, bounds for the Betti numbers of any $R$-module $M$ over the fiber product ring $R\times_TS$ are given  (Corollary \ref{cormainbetti}). In addition, we derive the Poincar\'e series and Betti numbers  of any $R$-module $M$ over the {\it amalgamated duplication ring} $R\bowtie I$. This special ring introduced by D'Anna \cite{danna}, improves  the notion of Nagata's idealization $R\ltimes I$ (also called trivial extension)(cf. \cite[p. 2]{nagata}) and  brought with it interesting geometric applications to curve singularities. Actually, D'Anna \cite{danna} showed that if $R$ is an algebroid curve with $h$ branches, then $R\bowtie I$ is also an algebroid curve with $2h$ branches, and more, a explicit form to construct Gorenstein algebroid curves was provided.
The characterization of the Betti numbers obtained is an important fact for the rest of the paper, especially for the structural results provided for  general fiber product rings.
 
 An interesting open problem in the theory of fiber product rings $R\times_TS$ is to give the description of their depth.    When  $T=k$, Lescot \cite{L81} has shown  the famous equality  $${\rm depth}(R\times_kS)= { \rm min}\{ {\rm depth}(R), {\rm depth}(S),1\},$$ but, if  $T\neq k$ is a $0$-dimensional ring, a similar formula is not expected (see Example \ref{exemplodepthfura}).This pathology brings us to the need to use other invariant to characterize the depth of $R\times_TS$. In this direction, Theorem \ref{depthgeral} provides some cases of the exact value of  ${\rm depth}(R\times_TS)$  using the grade of the extension in $R$, $S$ and $T$  of the maximal ideal  $\mathfrak{m}$ of $R\times_TS$. With arguments from local cohomology theory,  we generalize  the formula of Lescot \cite{L81}, as follows (see Corollary \ref{Lescotgeral}):

 \begin{cor}\label{Lescotgeral0} Let $R\times_TS$ be a fiber product ring. If  $\dim(T)=0$  and $\Gamma_{\mathfrak{m}R}(R) \subseteq {\rm ker}(\pi_R)$, then    
$${\rm depth}(R\times_TS)= { \rm min}\{ {\rm grade}(\mathfrak{m}R, R), {\rm grade}(\mathfrak{m}S, S), 1\}.$$
\end{cor}

Let us now briefly describe the contents of the paper. Section 2  establishes  our terminology and  some preparatory facts that will be needed throughout the paper. Section 3 deals of the Poincar\'e series any $R$-module $M$ over  $R\times_TS$ and their Betti numbers. We cite Theorem \ref{thm1} and  Corollary \ref{cormainbetti} as the main results of this section. Also, we define a new class of fiber product rings, called large, with contains any fiber product ring $R\times_kS$, for instance,  and we derive some results concerning their Poincar\'e series and Betti numbers. As a consequence, the Poincar\'e series of the amalgamated duplication ring \cite{danna} is given.     In Section 4,  we investigate the depth of $R\times_TS$ (Theorem \ref{depthgeral}) and we show some structural results of general fiber product and large fiber product rings (Theorem \ref{thmpregular1} and Proposition \ref{structlargefiber}). As the main consequence, we obtain that any fiber product ring is not a regular ring. Further, we provide some partial answers for the  Buchsbaum-Eisenbud-Horrocks and Total rank conjectures for modules with infinite projective dimension, using the  Betti numbers obtained in Section 3.

\section{Setup}

Let $(R,\mathfrak{m}_R,k)$, $(S,\mathfrak{m}_S,k)$ and $(T,\mathfrak{m}_T,k)$ be
 commutative local rings,  and let $R \stackrel{\pi_R}\twoheadrightarrow T \stackrel{\pi_S}\twoheadleftarrow S$ be
  surjective homomorphisms of rings.  The fiber product ring
$$
R \times_T S=\{(r,s)\in R\times S  \ \mid \ \pi_R(r)=\pi_S(s)  \},
$$
is a Noetherian local ring with maximal
ideal $\mathfrak{m}_R \times_{\mathfrak{m}_S}\mathfrak{m}_S$, residue field $k$, and it is a
subring of the usual direct product $R\times S$ (see \cite[Lemma 1.2]{AAM}). 





Let $\eta_R: R \times_T S\twoheadrightarrow R$
and $\eta_S:R \times_T S\twoheadrightarrow S$ be the natural
projections $(r,s)\mapsto r$ and $(r,s)\mapsto s$, respectively. Then $R \times_TS$ is represented as a pullback diagram:
\begin{equation}\label{dia}
\begin{CD}
R \times_T S & \stackrel{\eta_S}{\longrightarrow} & S\\
{\eta_R}{\downarrow} &  & {\downarrow}{\pi_S} \\
R & \stackrel{\pi_R}{\longrightarrow} & T.
\end{CD}
\end{equation}

\
\

 \begin{nota} \label{syzbnot}{\rm 
 Whenever we use the symbols $R$, $S$ or $T$, we  assume that $R$, $S$ and $T$ are non-zero Noetherian local rings with the maximal ideal $\m_R$, $\mathfrak{m}_S$ and $\mathfrak{m}_T$, respectively, and with $k$ being the common residue field.  
Every module over $R$, $S$ or $T$ is assumed to be finitely generated and non-zero.
The fiber product ring $R\times_TS$ is non-trivial, i.e.,   $R\neq T\neq S$, and their maximal ideal is denoted by $\m:= \m_R\times_{\m_T}\m_S$. Note that every $R$-module (or $S$-module) is an $R\times_T S$-module via Diagram \ref{dia}.

 For an $R$-module $M$,    let $\mu_R(M)$ denote the minimal number of generators of  $M$. Also, let $\edim (R)$ denote the {\it embedding dimension} of a local ring $R$, that is,  the minimal number of generators of the maximal ideal $\mathfrak{m}_R$, i.e.,   $\edim(R):=\dim_{k}(\mathfrak{m}_R/\mathfrak{m}_R^2)$. In general, one has  $\edim(R) \geq \dim (R)$ and the equality happens provided $R$ is  a regular  ring.

}
\end{nota}

Now, we  summarize some definitions and  known results for the rest of this work.

\begin{fact}\label{dimdepth}{ \rm For the fiber product ring  $R\times_TS$, \cite{AAM} and \cite{L81} (or \cite{paper}) provide
$$\dim (R\times_TS) = {\rm max} \{ \dim (R), \dim (S) \},$$
$${\rm depth} (R\times_TS) \geq {\rm min} \{{\rm depth} (R), {\rm depth} (S), {\rm depth}(T)+ 1\}.$$

If $T=k$, then
$${\rm depth} (R\times_kS) = {\rm min} \{{\rm depth} (R), {\rm depth} (S), 1\} .$$
 
 }
\end{fact}

\begin{defn}\label{defbasic} {\rm 
\begin{itemize}
\item[(i)]
Let $R$ be a  ring and let $M$ be an $R$-module. Consider the formal power series $$P_{M}^R(t):= \sum_{i\geq 0}\dim_{k} {\rm Tor}_i^R\left(M, k \right) t^i,$$
where $\dim_k(-)$ denotes the dimension of a vector space over $k$. The  series $P_{M}^R(t)$ is known as the {\it Poincaré series} of $M$ and the number $\beta_i^R(M):=\dim_{k} {\rm Tor}_i^R\left(M, k\right)$ is called the {\it $i$-th Betti number} of $M$.

\item[(ii)] (\cite{levin})  Let $f: R\to S$ be a surjective homomorphism. Then $f$ is said to be {\it large} provided, for any $S$-module $M$, considered as an $R$-module via $f$, the following equality happens
$$P_M^R=P_M^SP_S^R.$$

\item[(iii)] Let $F(t) = \sum
a_it^i$ and $G(t) = \sum
b_it^i$ be two power series in $t$. We say that $F(t)\succeq G(t)$ provided  $a_i \geq b_i$ for all $i$.

\end{itemize}
}    
\end{defn}

\begin{rem}\label{rem3.3} {\rm
\begin{itemize}
    \item[(i)] Let $M$
be a  $R$-module  with minimal number of generators $\mu(M)$. Then the equality
$$P^R_M(t)=\mu(M)+tP^R_{\Omega_1}(t)$$
holds, where $\Omega_1:=\Omega_1(M)$  denotes the first syzygy of $M$ (see \cite{DK75}).
    \item[(ii)] (\cite[3.11]{ramati})
A special case of fiber product rings $R\times_TS$ comes from the following diagram below
\begin{equation}\label{dia31}
\begin{tikzcd}
& R\arrow[dr, "\pi_R"]\\
& & T\\
S\arrow[r, "\pi"] & R\arrow[ru, "\pi_
R"]
\end{tikzcd}
\end{equation}
where $\pi$ is a surjective map. In this case,   the homomorphism $\eta_S: R\times_TS \to S$ is large, i.e., for any $S$-module $N$ 
$$P^{R\times_TS}_N(t)=P_N^R(t)P_S^{R\times_TS}(t).$$

\item[(iii)] If $R\times_TS$ is a fiber product ring, it is important to realize that  $\beta_1^S(T)\neq 0\neq \beta_1^S(T)$. In fact,  if $\beta_1^S(T)=0$ for instance, then $T$ is a free $S$-module. The surjective map $S\stackrel{\pi_S}\to T$  and the fact that $T\cong S^{\oplus r}$, implies that $r=1$, i.e., $S\cong T$. Therefore, $R\times_TS \cong R\times_TT $ is the trivial fiber product ring, and this is a contradiction  with Notation \ref{syzbnot}.
\end{itemize}
}
\end{rem}

\section{Poincar\'e series of general fiber product rings}

The main focus of this section is to show that the Poincar\'e series of any $R$-module $M$ over the fiber product ring $R\times_TS$ is bounded by a rational function (see Theorem \ref{thm1}). For this purpose, an important observation made by Levin \cite{levin} is given.

\begin{rem}\label{remlev} \cite[Proof of Theorem 1.1]{levin} {\rm Let $f: R\to S$ be a surjective homomorphism. The spectral sequence
${\rm Tor}_p^{S}(M,k) \otimes {\rm Tor}_q^{R}(S,k) \Rightarrow {\rm Tor}^{R}(M,k)$ provides
$$P_M^R\succeq P_M^SP_S^R.$$

}
\end{rem}

We are able to show the main result of this section. 

\begin{thm} \label{thm1} Let  $R\times_TS$ be a  fiber product ring. If $M$ is an $R$-module, then
$$P^{R\times_TS}_M(t) \succeq \frac{P^R_M(t)P^S_T(t)}
{P^R_T(t)+P^S_T(t)-P^R_T(t)P^S_T(t)}.$$
\end{thm}
\begin{proof}
From the exact sequence
\begin{equation}\label{th310}
    0\longrightarrow \ker\pi_S\longrightarrow R\times_TS\overset{\eta_ R}{\longrightarrow} R\longrightarrow 0,
\end{equation}
and Remark \ref {rem3.3} (i),  one has
\begin{equation}\label{th3200}
    P_R^{R\times_TS}(t)=1+tP^{R\times_TS }_{\ker{\pi_S}}(t).
\end{equation}
Since the map $R\times_TS\stackrel{\eta_S} \to S$ is surjective, Remark \ref{remlev} gives
\begin{equation}\label{th3300}
    tP^{R\times_ST }_{\ker{\pi_S}}(t) \succeq tP_{\ker{\pi_S}}^S(t)P_S^{R\times_TS}(t)=(P_T^S(t)-1)P^{R\times_T S}_S(t),
\end{equation}
where the last equality follows by the exact sequence
\begin{equation}
0\longrightarrow \ker\pi_S\longrightarrow S\overset{\pi_ R}{\longrightarrow} T\longrightarrow 0
    \end{equation}
and Remark \ref{rem3.3} (i).
Hence (\ref{th3200}) and (\ref{th3300}) provide  \begin{equation}\label{th340}
    P_R^{R\times_TS}(t)\succeq 1+(P_T^S(t)-1)P^{R\times _TS}_S(t).
\end{equation}

With a similar argument, the exact sequence 
$$0\longrightarrow \ker\pi_R\longrightarrow R\times_TS\overset{\eta_ R}{\longrightarrow} S\longrightarrow 0$$ and the fact the map  $R\times_TS\stackrel{\eta_R} \to R$ is surjective, yield
\begin{equation}\label{th3300}
 P_S^{R\times_TS}(t) \succeq  1+(P_T^R(t)-1)P^{R\times_T S}_R(t).
\end{equation}

Replacing   (\ref{th3300}) in   (\ref{th340}) one obtains
\begin{equation}\label{3600}
    P_R^{R\times_TS}(t)\succeq \frac{P_T^S(t)}{(P_T^R(t)+ P^S_T(t)- P^R_T(t)P^S_T(t))}.
\end{equation}
Since $P^{R\times_TS}_M(t)\succeq P^R_M(t)P^{R\times_TS}_R(t)$ for any $R$-module $M$, multiplying both sides of equation (\ref{3600}) by $P^R_M(t)$, we derive
$$P^{R\times_TS}_M(t) \succeq \frac{P^R_M(t)P^S_T(t)}
{P^R_T(t)+P^S_T(t)-P^R_T(t)P^S_T(t)}.$$
\end{proof}

Motivated by the definition given by Levin \cite{levin} (see Definition \ref{defbasic} (ii)),  below we introduce a new class of fiber product rings.

\begin{defn}\label{deflargefiber} {\rm  
 We say that the fiber product ring $R\times_TS$ is {\it large}  provided the maps $\eta_R$ and $\eta_S$ are both large (see Diagram \ref{dia}).  
}
\end{defn}

Note that the class of large fiber product rings contains, for instance, $R\times_kS$ (\cite[Proposition 3.1]{L81}) and  $R\times_TR$ (Remark \ref{rem3.3} (ii)). Also, if we assume that there is a  surjective ring homomorphism $R\to S$ and  ${\rm ker}(\pi_R)$ is a weak complete
intersection ideal in $R$, then $R\times_TS$ is also large  (\cite[Theorem 3.12]{ramati}).

As an immediate consequence of the proof of Theorem \ref{thm1} and definition of large fiber product rings, we derive the following result, that improves \cite{DK75} and \cite[Theorem A]{ramati}.
\begin{thm} \label{thm12} Let  $R\times_TS$ be a large  fiber product ring. If $M$ is an $R$-module, then
$$P^{R\times_TS}_M(t) = \frac{P^R_M(t)P^S_T(t)}
{P^R_T(t)+P^S_T(t)-P^R_T(t)P^S_T(t)}.$$
\end{thm}
\begin{rem}\label{remdefamalgament}\rm{
Note that, if we assume  that $I$ is a proper $R$-ideal, then $$R\bowtie I=\{(r, r+i)|\, r\in R,\, i\in I\}$$  \textit{is the amalgamated duplication of $R$ along $I$}, introduced by D'Anna \cite{danna}. In particular, if we assume  that $\pi_R=\pi_S:R\to R/I$ (Diagram \ref{dia}) are the canonical maps, then  $$R\times_{R/I}R=\{(r,s)|\, r+I=s+I\}= \{(r, r+i)|\, r\in R,\, i\in I\}= R\bowtie I.$$}  
\end{rem}

Since $R\times_{R/I}R= R\bowtie I$ is large, Theorem \ref{thm12} gives the  following consequence.

\begin{cor}\label{amalgament} Let  $R\bowtie I$ be the amalgamated duplication  ring of $R$ along to  the proper ideal $I$. If $M$ is an $R$-module, then  
$$P^{R\bowtie I}_{M}(t)=\frac{P^{R}_{M}(t)}
{2-P^{R}_{R/I}(t)}.$$
\end{cor}

Now, the focus is to give the bounds for the Betti numbers of any $R$-module over $R\times_TS$. Their proof is similar to \cite[Lemma 3.2]{paperbetti}, but we keep the proof for completeness.

\begin{cor}\label{cormainbetti} Let $R\times_TS$ be a fiber product ring and let $M$ be an  $R$-module. Then
\begin{itemize}
    \item[(i)] $\beta_0^{R\times_TS}(M) \geq  \beta_0^R(M).$ 
    \item[(ii)] 
    $\beta_1^{R\times_TS}(M)\geq  \beta_0^R(M)\beta_1^S(T)+\beta_1^R(M).$ In particular, $\edim(R\times_TS)\geq  \beta_1^S(T)+ \edim(R)$.
    \item[(iii)] 
    $\beta_2^{R\times_TS}(M)\geq \beta_0^R(M)\beta_1^R(T)\beta_1^S(T)+\beta_0^R(M)\beta_2^S(T)+\beta_1^R(M)\beta_1^S(T)+\beta_2^R(M).$
    
    \item[(iv)]  $\beta_i^{R\times_TS}(M)\geq \beta_i^R(M)$ for all $i\geq 0$.
    
\end{itemize}
\end{cor}
\begin{proof}
By Theorem \ref{thm1}, one has
$$P^{R\times_TS}_M(t)\succeq\frac{\left(\sum\limits_{i\geq 0} \beta_i^R(M) t^i\right)\left(\sum\limits_{i\geq 0} \beta_i^S(T) t^i\right)}
{\sum\limits_{i\geq 0} \beta_i^R(T) t^i+\sum\limits_{i\geq 0} \beta_i^S(T) t^i-\left(\sum\limits_{i\geq 0} \beta_i^R(T) t^i\right)\left(\sum\limits_{i\geq 0} \beta_i^S(T) t^i\right)}.$$

The equalities
$$\left(\sum\limits_{i\geq 0} \beta_i^R(M) t^i\right)\left(\sum\limits_{i\geq 0} \beta_i^S(T) t^i\right)=\sum\limits_{n\geq 0}\left(\sum\limits_{i= 0}^n \beta_i^R(M)\beta_{n-i}^S(T)\right)t^n$$ and
 
{\scriptsize{$$\sum\limits_{i\geq 0} \left(\beta_i^R(T) t^i+ \beta_i^S(T) t^i\right)-\sum\limits_{n\geq 0}\left(\sum\limits_{i= 0}^n \beta_i^R(T)\beta_{n-i}^S(T)\right)t^n=\sum\limits_{i\geq 0}\left(\beta_i^R(T)+\beta_i^S(T)-\sum\limits_{j= 0}^i \beta_j^R(T)\beta_{i-j}^S(T)\right)t^i,$$}}

give
$$P_M^{R\times_TS}(t)\succeq\frac{\sum\limits_{i\geq 0}\left(\sum\limits_{j= 0}^i \beta_j^R(M)\beta_{i-j}^S(T)\right)t^i}{\sum\limits_{i\geq 0}\left(\beta_i^R(T)+\beta_i^S(T)-\sum\limits_{j= 0}^i \beta_j^R(T)\beta_{i-j}^S(T)\right)t^i}.$$

The series  
$$\frac{1}{\sum\limits_{i\geq 0}\left(\beta_i^R(T)+\beta_i^S(T)-\sum\limits_{j= 0}^i \beta_j^R(T)\beta_{i-j}^S(T)\right)t^i}=\sum\limits_{i\geq 0}B_it^i,$$ is calculated as follows.  
Set $b_i:= \beta_i^R(T)+\beta_i^S(T)-\sum\limits_{j= 0}^i \beta_j^R(T)\beta_{i-j}^S(T)$. Note that $b_0=1$, and so    $B_0=1$. Since $b_0\neq 0$,    for all $i\geq 1$,  $B_i$ is given by the determinant

$$
B_i=\frac{1}{b_0^ii!}\begin{vmatrix}
0 & ib_1 & ib_2 &\dots& ib_i \\ 
0 & (i-1)b_0 & (i-1)b_1 &\dots& (i-1)b_{i-1} \\ 
0 & 0 & (i-2)b_0 &\dots& (i-2)b_{i-2} \\ 
\vdots & \vdots & \vdots &\ddots& \vdots\\
1 & 0 & 0 &\dots& 1
\end{vmatrix}.
$$
By induction on 
$B_n= \sum_{i=1}^n |b_i|B_{n-i}$, and  since $b_i\leq 0$ for $i\geq 1$,  note that $B_i \geq 0$.  Now, set 
$a_i=\sum\limits_{j= 0}^i \beta_j^R(M)\beta_{i-j}^S(T).$ Therefore
 $$P_M^{R\times_TS}(t)\succeq\left(\sum\limits_{i\geq 0}a_it^i\right)\left(\sum\limits_{i\geq 0}B_it^i\right)=\sum\limits_{n\geq 0}\left(\sum\limits_{i= 0}^n a_iB_{n-i}\right)t^n,$$ and this provides 
$\beta_{n}^{R\times_TS}(M)\succeq\sum\limits_{i= 0}^n a_iB_{n-i}.$ Now, the statements (i)-(iv) follows.
\end{proof}

As a biproduct of Theorem \ref{thm12} and a similiar proof of Corollary \ref{cormainbetti}, we derive the following: 

\begin{cor}\label{cormainbetti2} Let $R\times_TS$ be a large fiber product ring and let $M$ be an  $R$-module. Then
\begin{itemize}
    \item[(i)] $\beta_0^{R\times_TS}(M) = \beta_0^R(M).$ 
    \item[(ii)] 
    $\beta_1^{R\times_TS}(M)=  \beta_0^R(M)\beta_1^S(T)+\beta_1^R(M).$ In particular, $\edim(R\times_TS)=  \beta_1^S(T)+ \edim(R)$.
    \item[(iii)] 
    $\beta_2^{R\times_TS}(M)= \beta_0^R(M)\beta_1^R(T)\beta_1^S(T)+\beta_0^R(M)\beta_2^S(T)+\beta_1^R(M)\beta_1^S(T)+\beta_2^R(M).$

\end{itemize}
\end{cor}

\section{ Structural results and lower  bounds for Betti numbers}

A key ingredient for understanding the behavior of the structure of any ring, is its depth. In the investigation of fiber product rings $R\times_TS$,   to know the exact value of their depth is an important open problem. When $T=k$,  Lescot \cite{L81} has shown that  
\begin{equation}\label{depthlescot}
{\rm depth}(R\times_kS)= { \rm min}\{ {\rm depth}(R), {\rm depth}(S),1\},
\end{equation}
but in  general (see \cite[Lemma 1.5]{AAM}) only the following inequality is given, 
$${\rm depth}(R\times_TS)\geq  { \rm min}\{ {\rm depth}(R), {\rm depth}(S), {\rm depth}(T)+1\}. $$

The next example illustrates that an analogous formula to (\ref{depthlescot}) does not hold for any $0$-dimensional ring $T\neq k$. 
\begin{ex}\label{exemplodepthfura}{\rm (\cite[Example 4.2.14]{Ela})  Consider $R = k[x, y ]/(y^2)$, $S = k[x,y]/(x^2,xy)$ and $T = k[x,y]/
(x^2,xy, y^2)$. The fiber product ring $R\times_TS= k[x,y]/(x^2y)$ is  $1$-dimensional Gorenstein. Further, $\dim(T) = 0$,  $S$ is not Cohen-Macaulay (note that ${\rm depth}(S)=0$) and $R$ is $1$-dimensional Cohen-Macaulay ring. }
\end{ex}

Now, the main focus  is to give  the exact value of the depth of $R\times_TS$  in terms of the grade of certain ideals of $R$, $S$ and $T$. For this purpose, we first recall some main facts concerning the local cohomology theory, the key tool for the desired result. For details, see \cite{B-S} and \cite{pS}.  

Let  $M$ be a finitely generated module over the local ring  $(A,\mathfrak{m}_A,k)$. Let $H^i_{\mathfrak{a}}(M)$ denote the $i$-th local cohomology module of $M$ with respect to an $R$-ideal  $\mathfrak{a}$. The $H^i_{\mathfrak{a}}$ are the right derived functors of the
left exact $\mathfrak{a}$-torsion functor
\[
\Gamma_{\mathfrak{a}} (M) = \{x\in M \mid \mathfrak{a}^tx = 0 \text{ for some } t\in\mathbb{N}\}\,.
\]
 
\begin{enumerate}

\item Let $0\to M\to N\to L\to 0$ be a short exact sequence of finitely generated modules  over $A$.
Then there is an associated long exact sequence of cohomology modules
$$
\aligned  0\to & \Gamma_{\mathfrak{a}}(M)\to \Gamma_{\mathfrak{a}}(N)\to \Gamma_{\mathfrak{a}}(L)
\to H^1_{\mathfrak{a}}(M)\to H^1_{\mathfrak{a}}(N)\to H^1_{\mathfrak{a}}(L)\to\cdots
\\
& \cdots \to H^i_{\mathfrak{a}}(M)\to  H^i_{\mathfrak{a}}(N)\to H^i_{\mathfrak{a}}(L)\to H^{i+1}_{\mathfrak{a}}(M) \to
 \cdots\, .
\endaligned
$$

\item If $M$ is an $A$-module of finite length, then $\Gamma_{\m_A}(M) = M$, and $H^i_{\m_A}(M) = 0$ for every $i\ge 1$. Also, if $I \subseteq J$ are $A$-ideals, the $\Gamma_J(M) \subseteq \Gamma_I(M)$.

\item The length of the longest $M$-sequence contained in $\mathfrak{a}$ is  denoted by ${\rm grade}_R(\mathfrak{a},M)$. In the case that $\mathfrak{a}M\neq M$, then
$${\rm grade}(\mathfrak{a}, M) = \inf\{i \mid H^i_{\mathfrak{a}}(M) \ne 0\}.$$
\noindent When $\mathfrak{a}=\mathfrak{m}_A$, then the grade of $M$ in $\mathfrak{a}$ is called the ${\rm depth}$ of $M$.

\item   $H^i_{\mathfrak{a}}(M) = 0$ for every $i> \dim(M).$

\end{enumerate}

\begin{nota}{\rm For the rest of this paper, if $f:R\to R'$ is a homomorphism of Noetherian rings and $\mathfrak{a}$ is an ideal of $R$, set $\mathfrak{a}R'$ as the extension of $\mathfrak{a}$ to $R'$ under $f$.}
\end{nota}

We are able to show one of the main results of this section.

\begin{thm}\label{depthgeral} Let  $R\times_TS$ be a  fiber product ring.

\begin{itemize}
    \item[(i)] Suppose that ${\rm grade}(\mathfrak{m}T, T)=n \geq 0$,  ${\rm grade}(\mathfrak{m}R, R)> n$ and ${\rm grade}(\mathfrak{m}S, S)> n .$ Then ${\rm depth}(R\times_TS)= n+1$. 

    \item[(ii)] If ${\rm depth}(T)= 0$,  ${\rm grade}(\mathfrak{m}R, R)> 0$ and ${\rm grade}(\mathfrak{m}S, S)> 0,$ then ${\rm depth}(R\times_TS)= 1$.

    \item[(iii)] Assume that ${\rm grade}(\mathfrak{m}T, T)=n > 0$. If ${\rm depth}(R)=0$ or ${\rm depth}(S)=0,$ then ${\rm depth}(R\times_TS)= 0$.

\end{itemize}
\end{thm}
\begin{proof}
(i) Suppose that $n=0$. The short exact sequence (\cite[1.1.2]{AAM}) \begin{equation}\label{SeqexataLC}
 0\longrightarrow R \times_{T} S\longrightarrow R\oplus S\longrightarrow T \longrightarrow 0,   
\end{equation}
gives the
associated long exact sequence of local cohomology modules
\begin{equation}\label{longextLC}
0\to  \Gamma_{\m}(R\times_TS) \to \Gamma_{\m}(R) \oplus \Gamma_{\m}(S)  \to \Gamma_{\m}(T)
\to H^1_{\m}(R\times_TS)\to \cdots .
\end{equation}

The Independence Theorem of local cohomology modules (\cite[4.2.1]{B-S}) yields the isomorphisms $\Gamma_{\m}(R)\cong \Gamma_{\m R}(R)$, $\Gamma_{\m}(S)\cong \Gamma_{\m S}(S)$ and $\Gamma_{\m}(T)\cong \Gamma_{\m T}(T)$, and then from (\ref{longextLC}), the following long exact sequence
\begin{equation}\label{longextLC2}
0\to  \Gamma_{\m}(R\times_TS) \to \Gamma_{\m R}(R) \oplus \Gamma_{\m S}(S)  \to \Gamma_{\m T}(T)
\to H^1_{\m}(R\times_TS)\to \cdots .
\end{equation}

The hypothesis ${\rm grade}(\mathfrak{m}R, R)> 0$ and ${\rm grade}(\mathfrak{m}S, S)> 0$ provide $\Gamma_{\m_R}(R) = 0 = \Gamma_{\m_S}(S)$, and hence  $\Gamma_{\m}(R\times_TS) = 0$, from (\ref{longextLC2}).  In addition, one obtains that  $ \Gamma_{\m T}(T)
\to H^1_{\m}(R\times_TS)$ is an injective map. Since ${\rm grade}(\mathfrak{m}T, T)= 0$, one has  $\Gamma_{\m_T}(T)
\ne 0$. Therefore $H^1_{\m}(R\times_TS)\neq 0$, i.e., ${\rm depth}(R\times_TS)=1$, as desired.

Now, assume that ${\rm grade}(\mathfrak{m}T, T)=n >0$ and consider the long exact sequence 
\begin{equation}\label{longextLC3}
\cdots \to H^n_{\m}(R\times_TS) \to H^n_{\m R}(R)  \oplus \to H^n_{\m S}(S)  \to  H^n_{\m T}(T)
\to H^{n+1}_{\m}(R\times_TS)\to \cdots
\end{equation}
given in 
 (\ref{longextLC2}).
Since  ${\rm grade}(\mathfrak{m}R, R)> n$ and ${\rm grade}(\mathfrak{m}S, S)> n$, then $H^j_{\m R}(R)=0=H^j_{\m S}(S)$ for all $j\leq n$. Also, ${\rm grade}(\mathfrak{m}T, T)=n$ furnishes $H^i_{\m T}(T)=0$ for all $i<n$ and $H^n_{\m T}(T)\neq 0$. Therefore the long exact sequence (\ref{longextLC3}) provides $H^{l}_{\m}(R\times_TS)=0$ for all $l\leq n$ and that the map $H^n_{\m T}(T)
\to H^{n+1}_{\m}(R\times_TS)$ is injective. This yields 
$H^{n+1}_{\m}(R\times_TS) \neq 0$, and so ${\rm depth}(R\times_TS)=n+1$. 

(ii) By definition,  $0={\rm depth}(T)\geq {\rm grade}(\mathfrak{m}T, T)$. Therefore $ {\rm grade}(\mathfrak{m}T, T)=0$. Now, the result is a consequence of (i).

(iii) Since ${\rm grade}(\mathfrak{m}T, T)= n>0$, then $\Gamma_{\m T}(T)=0$ and so, (\ref{longextLC2}) gives the isomorphism
\begin{equation}\label{isoLCGamma}\Gamma_{\m}(R\times_TS) \cong \Gamma_{\m R}(R) \oplus \Gamma_{\m S}(S).\end{equation}
By definition, ${\rm depth}(R)\geq {\rm grade}(\mathfrak{m}R, R)$. So, if ${\rm depth}(R)=0$, then ${\rm grade}(\mathfrak{m}R, R)=0,$ and then $\Gamma_{\m R}(R)\neq 0$. Therefore (\ref{isoLCGamma}) gives $\Gamma_{\m}(R\times_TS) \neq 0$, i.e., ${\rm depth}(R\times_TS)=0$.  
\end{proof}

The next result improves  the famous formula given by Lescot \cite{L81}.

\begin{cor}\label{Lescotgeral} Let $R\times_TS$ be a fiber product ring. If  $\dim(T)=0$  and $\Gamma_{\mathfrak{m}R}(R) \subseteq {\rm ker}(\pi_R)$, then    
$${\rm depth}(R\times_TS)= { \rm min}\{ {\rm grade}(\mathfrak{m}R, R), {\rm grade}(\mathfrak{m}S, S), 1\}.$$
\end{cor}

\begin{proof} If we assume that  ${\rm grade}(\mathfrak{m}R, R)>0$ and ${\rm grade}(\mathfrak{m}S, S)>0$,  then Theorem \ref{depthgeral} (ii) yields ${\rm depth}(R
\times_kS)= 1$. Now, suppose that ${\rm grade}(\mathfrak{m}R, R)=0$ and ${\rm grade}(\mathfrak{m}S, S)>0$. Since $\dim T=0$, then  $\Gamma_{\mathfrak{m}_T}(T)=T$. Also, the the inclusion $\mathfrak{m}T \subseteq \mathfrak{m}_T$ provides $\Gamma_{\mathfrak{m}_T}(T)\subseteq \Gamma_{\mathfrak{m}T}(T)$ and then $\Gamma_{\mathfrak{m}T}(T)=T$.  The long exact sequence (\ref{longextLC2}) furnishes
\begin{equation}\label{longextLC00}
0\to  \Gamma_{\m}(R\times_TS) \to \Gamma_{\m_R}(R) \oplus \Gamma_{\m_S}(S)  \to T
\to H^1_{\m}(R\times_TS)\to \cdots .
\end{equation}
We want to show that $\Gamma_{\m}(R\times_TS)\neq 0$. Suppose that $\Gamma_{\m}(R\times_TS)=0$.  The assumptions ${\rm grade}(\mathfrak{m}R, R)=0$, ${\rm grade}(\mathfrak{m}S, S)>0$ imply $\Gamma_{\m_R}(R)\neq 0$ and $\Gamma_{\m_S}(S)= 0$, and then one has the injective map 
$0\to  \Gamma_{\m_R}(R)   \to T$, from (\ref{longextLC00}). Since $\Gamma_{\m_R}(R) \subseteq {\rm ker}(\pi_R)$ by hypothesis and $T\cong R/{\rm ker}(\pi_R)$, then $\Gamma_{\m_R}(R)=0$, a contradiction. Therefore $\Gamma_{\m}(R\times_TS)\neq 0$, i.e., ${\rm depth}(R\times_kT)=0$, and this proves the statement.  
\end{proof}

\begin{cor}{\rm (\cite{L81})}    ${\rm depth}(R\times_kS)= { \rm min}\{ {\rm depth}(R), {\rm depth}(S),1\}.$
\end{cor}
\begin{proof}
First, since that the  maximal ideal of  $R\times_kS$ is  $\mathfrak{m}= \mathfrak{m}_R\oplus \mathfrak{m}_S$ (see \cite{L81} or \cite{nasseh}), one has  $\mathfrak{m}R=\mathfrak{m}_R$ and $\mathfrak{m}S=\mathfrak{m}_S$. Hence,  ${\rm grade}(\mathfrak{m}R, R)= {\rm depth}(R)$ and ${\rm grade}(\mathfrak{m}S, S)= {\rm depth}(S)$. The definition of the fiber product ring $R\times_kS$ provides ${\rm ker}(\pi_R)=\m_R$. Therefore, the result follows by Corollary \ref{Lescotgeral}.
\end{proof}
A special case of fiber product rings are the amalgamated duplication  rings (Remark \ref{remdefamalgament}). 
When $R$ is Cohen-Macaulay and $I$ is a proper $R$-ideal, D'Anna \cite[Discussion 10]{danna} has shown that the depth of the amalgamated duplication ring  $R\bowtie I$ is equal of depth of $I$. However, if $R$ is not Cohen-Macaulay, the depth of  the amalgamated duplication ring is not known in the literature. In general, below we give a description of the depth of $R\bowtie I$. 
\begin{cor} Let  $I$ be a proper ideal of a ring $R$. 
\begin{itemize}
    \item[(i)] Suppose that ${\rm grade}(\mathfrak{m}R, R)>{\rm grade}(\mathfrak{m}R/I, R/I)=n \geq 0$. Then $${\rm depth}(R\bowtie I)= n+1.$$
    \item[(ii)] If
$\dim(R/I)=0$  and $\Gamma_{\mathfrak{m}R}(R) \subseteq I$, then    
$${\rm depth}(R\bowtie I)= { \rm min}\{ {\rm grade}(\mathfrak{m}R, R), 1\}\leq 1.$$ 
\end{itemize}

\end{cor}
\begin{proof}  Since $R\bowtie I=R\times_{R/I}R$ (Remark \ref{remdefamalgament}), the proof of  (i) and (ii) follows by Theorem \ref{depthgeral} (i) and Corollary \ref{Lescotgeral}, respectively.
\end{proof}

\subsection*{Structural Results} In this subsection we focus on the structural results of general fiber product rings, using the characterization of Betti numbers and depth obtained. Recall  that  a ring $R$ is a {\it hypersurface} provided there is a presentation $\widehat{R}= Q/I$, where $\widehat{R}$ is the $\mathfrak{m}$-adic completion of $R$, $Q$ is a regular local ring and $I$ is a principal ideal or, equivalently, if $\edim(R)- {\rm depth} (R)\leq 1$. The next result can be compared with \cite[Proposition 1.7]{AAM}, \cite[Proposition 4.19]{shirogoto} and \cite[Remark 3.1]{grow} (or \cite[Corollary 2.7]{nstv}).

\begin{thm}\label{thmpregular1} Let $R\times_TS$ be a fiber product ring.
\begin{itemize}
    \item[(i)]   $R\times_TS$  is not a regular ring.
    
    \item[(ii)]  If $R\times_T S$ is a hypersurface, then $R$ is regular and $\beta_1^S(T)=1$.

    \item[(iii)]  Suppose  that $\dim(R\times_TS)>0$ and ${\rm grade}(\mathfrak{m}T, T)=n > 0$. If ${\rm depth}(R)=0$ or ${\rm depth}(S)=0,$  then $R\times_TS$ is not Cohen-Macaulay.

\item[(iv)]  Assume that ${\rm grade}(\mathfrak{m}T, T)=n \geq 0$,  ${\rm grade}(\mathfrak{m}R, R)> n$ and ${\rm grade}(\mathfrak{m}S, S)> n .$ Then  $R\times_TS$ is Cohen-Macaulay if and only if     $R$ and $S$ are Cohen-Macaulay rings satisfying $\dim(R)=\dim(S)=\dim(R\times_TS)=n+1$, and ${\rm grade}(\mathfrak{m}T, T)= \dim(R)-1.$
 
\end{itemize}
\end{thm}
\begin{proof}  Without loss of generality, by Fact \ref{dimdepth} we  may assume  that $\dim (R\times_TS) =\dim (R)$.

(i) Suppose that $R\times_TS$ is regular.  Then, by Corollary \ref{cormainbetti} (ii), 
$$\dim (R)=\dim (R\times_TS)= \edim(R\times_TS) \geq \beta_1^S(T)+\edim(R).$$  Then $\beta_1^S(T)=0$ because  $\edim(R)\geq \dim (R)$. This gives a contradiction by Remark \ref{rem3.3} (iii).

(ii) Since $R\times_TS$ is a hypersurface, then ${\rm depth} (R\times_TS)= \dim(R\times_TS)=\dim (R)$ and $\edim(R\times_TS)-{\rm depth} (R\times_TS)\leq 1$. By Corollary  \ref{cormainbetti} (ii), one has
$$\beta_1^S(T)+\edim(R)-{\rm dim} (R)\leq \edim(R\times_TS)-{\rm depth} (R\times_TS)\leq 1.$$
The facts   $\beta_1^S(T)\neq 0$ (Remark \ref{rem3.3} (iii)) and $\edim (R)\geq \dim (R)$ provide   $\beta_1^S(T)=1$ and $\edim (R)={\rm dim}(R)$ (i.e. $R$ is regular), as desired.

(iii) By hypothesis, Theorem \ref{depthgeral} (iii) provides ${\rm depth}(R\times_TS)=0$. Since $R\times_TS$ has positive dimension, the desired statement follows.

(iv)  Since ${\rm depth}(R\times_TS)=n+1$  (Theorem \ref{depthgeral} (i)), $R\times_TS$ is Cohen-Macaulay if and only if  $\dim(R\times_TS)= n+1$. This gives the converse of the desired statement.  

By the fact ${\rm grade}(\mathfrak{m}R, R)> {\rm grade}(\mathfrak{m}T, T)$, there is a positive integer $\alpha$ such that 
\begin{equation}\label{igualdade1}{\rm grade}(\mathfrak{m}R, R)- \alpha={\rm grade}(\mathfrak{m}T, T).
\end{equation}
In addition, Theorem \ref{depthgeral} (i) gives 
${\rm depth}(R\times_TS)= {\rm grade}(\mathfrak{m}T, T)+1$. Hence, 
\begin{equation}\label{depthalpha}{\rm depth}(R\times_TS)= {\rm grade}(\mathfrak{m}R, R) -\alpha+1.
\end{equation}

Now, (\ref{depthalpha}) and the fact that $R\times_TS$ is Cohen-Macaulay provide \begin{equation}\label{depthalpha2}\dim(R)=\dim(R\times_TS)= {\rm depth}(R\times_TS)= {\rm grade}(\mathfrak{m}R, R) -\alpha+1.
\end{equation}
But ${\rm grade}(\mathfrak{m}R, R)-\alpha+1=\dim (R)\geq {\rm depth}(R)$ and ${\rm depth}(R) \geq  {\rm grade}(\mathfrak{m}R, R)$. Then, $-\alpha+1=0$, and so $\alpha =1$. This yields ${\rm depth}(R) =  {\rm grade}(\mathfrak{m}R, R)$, and hence from (\ref{depthalpha2}) one obtains ${\rm depth}(R)=\dim(R)$, i.e., $R$ is Cohen-Macaulay. Furthermore, from (\ref{igualdade1}) one has ${\rm grade}(\mathfrak{m}T, T)= \dim(R)-1$.

Similarly to (\ref{depthalpha}) and (\ref{depthalpha2}), we derive ${\rm depth}(R\times_TS)= {\rm grade}(\mathfrak{m}S, S)- \beta+1$, and
$$\dim(R)=\dim(R\times_TS)= {\rm depth}(R\times_TS)= {\rm grade}(\mathfrak{m}S, S)- \beta+1.$$

Since $R$ is Cohen-Macaulay and $\dim(R\times_TS)=
\dim(R)$, then ${\rm grade}(\mathfrak{m}S, S)- \beta+1=\dim (R)= {\rm depth}(R) \geq \dim(S) \geq {\rm depth}(S)$. By definition,  ${\rm depth}(S) \geq  {\rm grade}(\mathfrak{m}S, S)$, and this gives $-\beta+1=0$. This furnishes $\beta=1$, and since the dimension of $R\times_TS$ is the maximum of $\dim(R)$ and $\dim(S)$ (Fact \ref{dimdepth}), we get  
$$\dim(R\times_TS)=\dim(R)= {\rm depth}(R)= {\rm depth}(R)=\dim(S).$$ Therefore $S$ is also Cohen-Macaulay satisfying $\dim(S)=\dim(R)= \dim(R\times_TS)$.
\end{proof}

\begin{cor}\label{cor4.7}  Suppose that ${\rm depth}(T)= 0$,  ${\rm grade}(\mathfrak{m}R, R)> 0$ and ${\rm grade}(\mathfrak{m}S, S)> 0.$ Then $R\times_TS$ is Cohen-Macaulay if and only $R$ and $S$ are  Cohen-Macaulay rings satisfying  $\dim(R)=\dim(S)=\dim(R\times_TS)= 1$.
\end{cor}
\begin{proof}
Since  $0={\rm depth}(T)\geq {\rm grade}(\mathfrak{m}T, T)$, one has ${\rm grade}(\mathfrak{m}T, T)=0$. Now, the result follows by Theorem \ref{thmpregular1} (iv).    
\end{proof}

\begin{ex}{\rm  Consider the Cohen-Macaulay fiber product  ring $R\times_TS= k[x,y]/(x^2y)$ given in  Example  \ref{exemplodepthfura}. Note that $T$ is $0$-dimensional,  ${\rm depth}(S)=0$ (and then ${\rm grade}(\mathfrak{m}S,S)=0$).  Also, $S$ is not Cohen-Macaulay ring. This example illustrates that the conditions of the previous result can not be improved. 
}
\end{ex}


When the fiber product ring $R\times_TS$ is large, other structural results are provided. The next result improves  \cite[Propositions 3.21, 4.7, 4.11 and Corollary 4.9]{paperbetti}, because if $R\times_TS$ satisfies $(\ast)$ (\cite{paperbetti}) then the fiber product is large. Their proof is similar and we omit it here. Recall that a ring $R$ is  called a {\it complete intersection} if its completion in the $\mathfrak{m}$-adic topology is the quotient of a regular ring $Q$ by a $Q$-regular sequence $\underline{x}=x_1,\ldots,x_c$ contained in the maximal ideal of $Q$. 

\begin{prop}\label{structlargefiber} Let $R\times_TS$ be a large fiber product ring. 
\begin{itemize}
\item[(i)] Suppose that $R$ is  a complete intersection.  Then $R\times_T S$ is a complete intersection provided 
$$\frac{\beta_1^S(T)^2+\beta_1^S(T)}{\beta_1^R(T)\beta_1^S(T)+\beta_2^S(T)}=2.$$

\item[(ii)] Assume  that $R\times_TS$ is Cohen-Macaulay.  Then $R\times_T S$ is a hypersurface if and only if $R$ is regular and $\beta_1^S(T)=1$.

\item[(iii)] If $R\times_TS$ is a Cohen-Macaulay  ring and  $\beta_{2}^{R\times_TS}(k) \leq  \beta^{R\times_TS}_{1}(k)$, then 
    $R\times_T S$ is  an hypersurface and $\pd_S (T)=1$.

\item[(iv)] Suppose that  $R\times_TS$ is Cohen-Macaulay. If $\beta_1^{R\times_TS}(M) > \beta^{R\times_TS}_0(M)$ for any $R$-module $M$, then  $R\times_T S$ is not a hypersurface or $M$ is not a free $R$-module.

\end{itemize}

\end{prop}

Regarding the structure of the amalgamated duplication ring $R\bowtie I$, few results are known in the literature. In fact, it is only characterized when $R\bowtie I$ is Cohen-Macaulay or Gorenstein (see \cite[Theorem 1.8 and 1.9 - 1.9.4]{AAM} for details). Our previous results provide a refinement to this investigation, showing that $R\bowtie I$ is always non-regular, as well as when it is a hypersurface or a complete intersection ring.

\begin{cor} Let  $R\bowtie I$ be the amalgamated duplication of $R$ along $I$.

\begin{itemize}
    \item[(i)] $R\bowtie I$ is not a regular ring.
    
    \item[(ii)]  Suppose that $R\bowtie I$ is Cohen-Macaulay. Then $R\bowtie I$ is a hypersurface if and only if $R$ is regular and $\beta_1^R(R/I)=1$.
    
    \item[(iii)] Suppose that   $R$ is  a complete intersection. Then $R\bowtie I$ is a complete intersection if $\beta_2^R(R/I)=0$ and  $\beta_1^R(R/I)=1$.

\item[(iv)] Assume that $\dim(R/I)=0$ and ${\rm grade}(\mathfrak{m}R, R)> 0$. Then $R\bowtie I$ is Cohen-Macaulay if and only if $R$ is $1$-dimensional Cohen-Macaulay.
\end{itemize}
\end{cor}
\begin{proof} 
Since $R\bowtie I= R\times_{R/I}R$ (Remark \ref{remdefamalgament}), (i)-(iii) follows by  Theorem \ref{thmpregular1} and Proposition \ref{structlargefiber}. For (iv), note that  Fact \ref{dimdepth} provides that $\dim(R)=\dim(R\bowtie I)$. Now, the desired statement is a consequence of  Corollary \ref{cor4.7}.    
\end{proof}

\subsection*{Applications: Lower bounds of Betti numbers}
The last part of this paper is devoted to give some lower bounds for the Betti numbers of any $R$-module $M$ over the  fiber product $R\times_TS$. This investigation is motivated by the famous Buchsbaum-Eisenbud-Horrocks (\cite{BE, Ho}) and Total Rank Conjectures (\cite{AB}) that state the following:

\begin{conjec} [Buchsbaum-Eisenbud-Horrocks] Let $(R, \mathfrak{m}, k)$ be a
$d$-dimensional Noetherian local ring, and let $M$ be a finitely generated nonzero  $R$-module. If $M$ has finite length and finite projective dimension, then for each $i\geq 0$ the $i$-th Betti number of $M$ over $R$ satisfies the inequality
$$\beta_i^{R}(M)\geq \binom{d}{i}.$$
\end{conjec} 

This conjecture, from now on referred to as the (BEH) Conjecture, has a positive answer for local rings with dimension $\leq 4$ (see \cite{AB}), but for  larger dimensions the problem is still open. Serre (\cite[Lemma 9]{HU}), in a famous result in the literature, showed that this conjecture is true when $R$ is a  regular ring and $M=k$. Some other positive  answers in certain cases are provided, for instance, in \cite{Chang, Cha, EG}.
In \cite{AB}, Avramov and Buchweitz introduced a weaker version of the (BEH) Conjecture,  that asserts the following:

\begin{conjec} [Total Rank Conjecture] Let $(R, \mathfrak{m}, k)$ be a
$d$-dimensional Noetherian local ring, and let $M$ be a finitely generated nonzero  $R$-module. If $M$ has finite length and finite projective dimension, then 
$$\sum_{i\geq 0}\beta_i^{R}(M)\geq 2^d.$$
\end{conjec}

To simplify things, the Buchsbaum-Eisenbud-Horrocks and Total rank conjectures will be 
 now on referred to as the (BEH) and (TR) conjectures,  respectively. Note that both conjectures can be considered for modules with infinite projective dimension. In this direction,   Tate  \cite{tate}
has shown that  $\beta_i^R(k)=2^d$ for all $i\geq \dim R$, when $R$ is a hypersurface.

For the rest of this section, we say that  the ring $R$ satisfies the (BEH) Conjecture (respectively (TR) Conjecture), if for all $R$-module $M$ of finite length and finite projective dimension over $R$, we have $\beta_i^{R}(M)\geq \displaystyle\binom{d}{i}$ (respectively $\sum_{i\geq 0}\beta_i^{R}(M)\geq 2^d)$.

 The next result improves \cite[Theorem 1.3 and Proposition 3.17]{paperbetti}.
\begin{thm}\label{BEHTR} Let $R\times_TS$ be a $d$-dimensional ring, and let $M$ be an $R$-module. 
\begin{enumerate}

\item[(i)] If $R$ satisfies the {\rm (BEH)} Conjecture and $M$ has finite length and finite projective dimension over $R$, then for all $1\leq i \leq {\rm depth} (R)$
    $$\beta_i^{R\times_TS}(M)\geq \binom{d}{i}.$$

    \item[(ii)] If $R$ satisfies the {\rm (TR)} Conjecture and $M$ has finite length and finite projective dimension, then  $$\sum_{i= 0}^{ {\rm depth}(R)}\beta_i^{R\times_TS}(M)\geq 2^{d}.$$
    
\end{enumerate}
\end{thm}
\begin{proof} First,  we may assume that $\dim (R\times_TS) = \dim (R)$ (Fact \ref{dimdepth}).  By Corollary \ref{cormainbetti} (iv), one has
  $\beta_i^{R\times_TS}(M)\geq\beta_i^R(M).$ Therefore, if $R$ satisfies the (BEH) Conjecture, one obtains 
$$\beta_i^{R\times_TS}(M)\geq \beta_i^R(M)\geq \binom{d}{i},$$ for all $1\leq i\leq {\rm depth} (R)$. The proof of (ii) is similar to (i).
\end{proof}

It is important to realize that in the previous result, $M$ has infinite projective dimension over $R\times_TS$. In fact, since $\beta_1^S(T)\neq 0\neq \beta_1^S(T)$ by Remark \ref{rem3.3} (iii), all the Betti numbers given in Corollary \ref{cormainbetti} are greater than zero. So, Theorem \ref{BEHTR} furnishes a positive answer for the (BEH) and (TR) conjectures for modules that do not have finite projective dimension.

Further, since the (BEH)  conjecture is true  for the residue field $k$ over a regular ring  $R$, and any ring that satisfies the (BEH) conjecture also satisfies the (TR) conjecture,  we derive the following consequence.
\begin{cor} Let $R\times_TS$ be a $d$-dimensional fiber product ring. 
\begin{enumerate}

\item[(i)] If $R$ is regular, then for all $1\leq i \leq {\rm depth} (R)$
    $$\beta_i^{R\times_TS}(k)\geq \binom{d}{i}.$$

    \item[(ii)] If $R$ is regular, then  $$\sum_{i= 0}^{ {\rm depth}(R)}\beta_i^{R\times_TS}(k)\geq 2^{d}.$$
    
\end{enumerate}
\end{cor}

Another interesting consequence of Theorem \ref{BEHTR}  is due to a   result of Walker (\cite[Theorem 2 (1)]{walker}) which shows that (TR) Conjecture is true for  complete intersection rings where the characteristic  of the residue field $k$ is not two. It is important to realize that if $R$ is a complete intersection, the fiber product ring $R\times_TS$ is not always a complete intersection (see Theorem \ref{thmpregular1} (iv) and Proposition \ref{structlargefiber}). Therefore,   the next result  furnishes a new class of rings that satisfies the (TR) conjecture.    

\begin{cor}\label{walker0} Let $R\times_TS$ be a $d$-dimensional fiber product ring.  
Suppose that $R$ is the quotient of a regular local ring by a regular sequence of elements and the  characteristic of $k$ is different from 2. If $M$ is an $R$-module with finite projective dimension and finite length over $R$, then 
$$\sum_{i=0}^d\beta_i^{R\times_TS}(M)\geq 2^d.$$
\end{cor}

Recently,  VandeBogert and  Walker (\cite[Theorem 1.7]{New}) showed the true  of the Total Rank Conjecture for rings of characteristic two. As a consequence of this result and Theorem \ref{BEHTR}, we derive the following: 

\begin{cor}
   Let $R\times_TS$ be a $d$-dimensional fiber product ring.  Suppose that $R$ has characteristic $2$. If $M$ is an $R$-module with finite projective dimension and finite length over $R$, then 
    $$\sum_{i=0}^d\beta_i^{R\times_TS}(M)\geq 2^d.$$
\end{cor}

\section{Acknowledgments}
The authors would like to thank Victor Hugo Jorge Pérez for the kind suggestions  
 and  the encouragement  over the years. Also, the authors also wish to express their gratitude to the anonymous referee for their contributions to improving the paper.

\end{document}